\documentclass[11pt]{scrartcl}

\usepackage{amsmath,amsthm,amssymb,amsfonts}
\usepackage[colorlinks=true,urlcolor=blue,citecolor=blue,linkcolor=blue]{hyperref}

\usepackage{array}

\def \Uad{{U_{\text{ad}}}}

\numberwithin{equation}{section}

\newcolumntype{x}[1]{>{\centering\hspace{0pt}\arraybackslash}p{#1}}

\makeatletter
\renewcommand*{\env@matrix}[1][*\c@MaxMatrixCols c]{%
  \hskip -\arraycolsep
  \let\@ifnextchar\new@ifnextchar
  \array{#1}}
\makeatother

\begin{document}

\title{A sharp regularization error estimate for bang-bang solutions for an iterative Bregman regularization method for optimal control problems}

\author{
Frank P\"orner\footnote{Department of Mathematics,  University of W\"urzburg, Emil-Fischer-Str. 40, 97074 W\"urzburg, Germany, E-mail: frank.poerner@mathematik.uni-wuerzburg.de}
}
\date{}
\maketitle

\begin{abstract}
  In the present work, we present numerical results for an iterative method for solving an optimal control problem with inequality contraints. The method is based on generalized Bregman distances. Under a combination of a source condition and a regularity condition on the active sets convergence results are presented. Furthermore we show by numerical examples that the provided a-priori estimate is sharp in the bang-bang case.
\end{abstract}

\section{Introduction}
In this article we consider optimization problems of the following form
\begin{equation}\label{eq:main_problem}\tag{\textbf{P}}
    \text{Minimize} \quad \frac{1}{2}\|Su - z\|_Y^2  \quad \text{such that} \quad u_a \leq u \leq u_b \quad \text{a.e. in } \Omega
\end{equation}
which can be interpreted both as an optimal control problem or as an inverse problem. Here $\Omega \subseteq \mathcal{R}^n$, $n \geq 1$ is a bounded, measurable set, $Y$ a Hilbert space, $z \in Y$ a given function. The operator $S : L^2(\Omega) \to Y$ is linear and continuous. Here, the interesting  situation is, when $z$ cannot be reached due to the presence of the control constraints (non-attainability). The set of admissible functions is abbreviated by $\Uad := \{u \in L^2(\Omega): \, u_a \leq u \leq u_b\}$. We are interested in an iterative method to solve \eqref{eq:main_problem} based on generalized Bregman distances. In \cite{wachsmuth2016} the algorithm was analysed under a suitable regularity assumption. Here we recall the most important results, followed by numerical results.

\section{Bregman iteration}
The Bregman distance associated with the regularization functional $J: u \mapsto \frac{1}{2}\|u\|_{L^2(\Omega)}^2 + I_\Uad(u)$ is defined as $D^\lambda(u,v) := J(u) - J(v) - (u-v, \lambda)$ where $\lambda \in \partial J(v)$. In the following $(\alpha_k)_k$ denotes a positive, bounded sequence of real numbers. The algorithm is given by:\\

\begin{tabular}{@{}|cl@{}}
\multicolumn{2}{|l}{Let $u_0 = P_\Uad(0)  \in \Uad$, $\lambda_0 = 0  \in \partial J(u_0)$ and $k=1$.}\\
1. & Solve for $u_k$:
\(\displaystyle \quad \text{Minimize} \quad \frac{1}{2}\|Su-z\|_Y^2 + \alpha_{k}  D^{\lambda_{k-1}}(u,u_{k-1}).\)\\
2. & Set $\lambda_k := \sum\limits_{i=1}^k \frac{1}{\alpha_i} S^\ast(z-Su_i) \in \partial J(u_k)$.\\
3. & Set $k:=k+1$, go back to 1.
\end{tabular}\\

The algorithm is well-defined due to the convexity of $D^\lambda(\cdot, v)$ with respect to the first argument (see \cite{wachsmuth2016} and the references therein).

\section{A-priori error estimates}
Let $u^\dagger$ be a solution of \eqref{eq:main_problem} and $p^\dagger=S^\ast( Su-z)$ be the adjoint state, then $(p^\dagger, u - u^\dagger) \geq 0, \, \forall u \in \Uad$ is satisfied. To derive our error estimates furthermore assume that there exists a set $I \subset \Omega$, $w \in Y$ and $\kappa, c > 0$ such that $I \supset \{ x \in \Omega: \; p^\dagger(x) = 0 \}$ holds. In addition assume that $\chi_I u^\dagger = \chi_I P_\Uad (S^\ast w)$ and $S^\ast w \in L^\infty(\Omega)$ holds. On the set $A := \Omega \setminus I$ we assume that the following structural assumptions $|\{ x\in A: \; 0 < |p^\dagger(x)| < \varepsilon  \}| \leq c \varepsilon^\kappa \,\,\, \forall \varepsilon > 0$ holds.

Under this regularity assumption strong convergence of the iterates $(u_k)_k$ can be established together with the a-priori error estimate
$$\|u^\dagger - u_k\|_{L^2(\Omega)}^2 \leq \mathcal{O} \left( \gamma_k^{-1} + \gamma_k^{-1} \sum\limits_{j=1}^k \alpha_j^{-1} \gamma_j^{- \kappa} \right),$$
with the abbreviation $\gamma_k := \sum_{j=1}^k \alpha_j^{-1}$. For details - both for the regularity assumption and the convergence - we refer to \cite{wachsmuth2016}. For the special choice of a constant sequence $\alpha_k = \alpha > 0$ and $\kappa < 1$ the a-priori estimate reduces to $\|u^\dagger - u_k\|^2 = \mathcal{O} \left( k^{-\kappa} \right)$ and to $\|u^\dagger - u_k\|^2 = \mathcal{O} \left( k^{-1} \log(k) \right)$ for $\kappa = 1$.

\section{Numerical examples}
In this section we present numerical results. The implementation is done in FEniCS \cite{fenics2016} with a semi-smooth Newton solver (see \cite{beuchler2010}). We use constant $\alpha_k = \alpha$ and compute the numerical approximation
$$\kappa_k :=  \frac{1}{\log(2)}\log \left( \frac{\|u_{k/2} - u^\dagger\|_{L^2(\Omega)}^2}{\|u_k-u^\dagger\|_{L^2(\Omega)}^2}  \right)$$
for bang-bang test examples ($A = \Omega$). Here our operator $y=Su$ is chosen to be the solution of the equation $-\Delta y = u$ in $\Omega$ and $y = 0$ on $\partial \Omega$. First we compute 1D examples with $\kappa = 1$, $\kappa = \frac{1}{2}$, and $\kappa =\frac{1}{3}$ for different mesh sizes $h$. The results are listed in Table 1,2 and 3 respectively. For the details of the construction of bang-bang examples with given adjoint state $p^\dagger$ we refer to \cite[Chapter 2.9]{troeltzsch2010}. To obtain $\kappa = 1$ we use $p^\dagger(x)=\sin(\pi x)$ on $\Omega=[-1,1]$. The other examples can be constructed using polynomials and limiting the slope near the zeros. For $\kappa = \frac{1}{3}$ we use $p^\dagger(x) = x(1-x)(3x-1)^3$ on $\Omega=[0,1]$.

Second we present a 2D bang-bang example, namely $p^\dagger(x,y) = \sin(2 \pi x) \sin(2 \pi y)$ on $\Omega = [0,1]^2$. Numerical estimates indicate  $\kappa = 1$, which is supported by our numerical results. Note that if the grid is too coarse the discretization error is dominating the regularization error, leading to unreliable results for $\kappa_k$. In all cases we obtain $\kappa_k \approx \kappa$ for $k$ large and $h$ small enough, indicating that our a-priori error estimate is sharp for the bang-bang case.

\vspace{0.75cm}

\centering
\begin{minipage}{70mm}

\begin{tabular}{@{} c | x{1.01cm} x{1.01cm} x{1.01cm} x{1.01cm} x{1.01cm} @{}}
$h$ & $10^{-3}$ & $10^{-4}$ & $10^{-5}$ & $10^{-6}$ \\\hline
$k$ & $\kappa_k$ & $\kappa_k$ & $\kappa_k$ & $\kappa_k$\\
\hline
4    & 0.646  & 0.602  & 0.601 & 0.601\\
8    & 0.839  & 0.752  & 0.750 & 0.750\\
16   & 1.027  & 0.860  & 0.856 & 0.857\\
32   & 1.211  & 0.927  & 0.922 & 0.923\\
64   & 1.229  & 0.960  & 0.958 & 0.960\\
128  & -0.001 & 0.945  & 0.975 & 0.979\\
256  & -0.004 & 0.786  & 0.980 & 0.989\\
512  & -0.020 & 0.271  & 0.972 & 0.991\\
1024 & -0.081 & -0.054 & 0.938 & 0.978\\
2048 & -0.217 & -0.149 & 0.826 & 0.919\\
\hline
\end{tabular}
\captionof{table}{1D example 1 ($\kappa = 1$).}

\end{minipage}
\hfil
\begin{minipage}{70mm}

\begin{tabular}{@{} c | x{1.01cm} x{1.01cm} x{1.01cm} x{1.01cm} x{1.01cm} @{}}
$h$ & $10^{-3}$ & $10^{-4}$ & $10^{-5}$ & $10^{-6}$ \\\hline
$k$ & $\kappa_k$ & $\kappa_k$ & $\kappa_k$ & $\kappa_k$\\
\hline
4    & 0.522 & 0.520 & 0.520 & 0.520\\
8    & 0.648 & 0.644 & 0.643 & 0.643\\
16   & 0.641 & 0.635 & 0.634 & 0.634\\
32   & 0.646 & 0.636 & 0.635 & 0.635\\
64   & 0.639 & 0.624 & 0.622 & 0.622\\
128  & 0.625 & 0.605 & 0.602 & 0.602\\
256  & 0.609 & 0.585 & 0.581 & 0.581\\
512  & 0.591 & 0.567 & 0.562 & 0.562\\
1024 & 0.571 & 0.553 & 0.547 & 0.546\\
2048 & 0.545 & 0.542 & 0.534 & 0.534\\
\hline
\end{tabular}
\captionof{table}{1D example 2 ($\kappa = \frac{1}{2}$).}

\end{minipage}

\vspace{0.75cm}

\centering
\begin{minipage}{70mm}

\begin{tabular}{@{} c | x{1.01cm} x{1.01cm} x{1.01cm} x{1.01cm} x{1.01cm} @{}}
$h$ & $10^{-3}$ & $10^{-4}$ & $10^{-5}$ &  $10^{-6}$  \\\hline
$k$ & $\kappa_k$ & $\kappa_k$ & $\kappa_k$ & $\kappa_k$\\\hline
4    & 0.286 & 0.286 & 0.286 & 0.286\\
8    & 0.312 & 0.312 & 0.312 & 0.312\\
16   & 0.325 & 0.327 & 0.328 & 0.328\\
32   & 0.329 & 0.337 & 0.338 & 0.338\\
64   & 0.321 & 0.339 & 0.340 & 0.341\\
128  & 0.301 & 0.338 & 0.340 & 0.340\\
256  & 0.272 & 0.335 & 0.338 & 0.339\\
512  & 0.236 & 0.332 & 0.337 & 0.338\\
1024 & 0.193 & 0.328 & 0.336 & 0.337\\
2048 & 0.132 & 0.322 & 0.335 & 0.336\\
\hline
\end{tabular}
\captionof{table}{1D example 3 ($\kappa = \frac{1}{3}$).}

\end{minipage}
\hfil
\begin{minipage}{70mm}

\begin{tabular}{@{} c | x{1.01cm} x{1.01cm} x{1.01cm} x{1.01cm} x{1.01cm} @{}}
$DOF$ & $10^{4}$ & $10^{5}$ & $10^{6}$ & $2 \cdot 10^{6}$ \\\hline
$k$ & $\kappa_k$ & $\kappa_k$ & $\kappa_k$ & $\kappa_k$\\
\hline
4    & 0.509  & 0.472 & 0.458 & 0.456\\
8    & 0.676  & 0.622 & 0.595 & 0.592\\
16   & 0.789  & 0.759 & 0.711 & 0.705\\
32   & 0.720  & 0.885 & 0.803 & 0.791\\
64   & 0.411  & 1.000 & 0.884 & 0.863\\
128  & 0.216  & 1.027 & 0.968 & 0.935\\
256  & 0.145  & 0.855 & 1.039 & 1.012\\
512  & 0.166  & 0.556 & 1.011 & 1.039\\
1024 & 0.129  & 0.295 & 0.805 & 0.936\\
2048 & -0.045 & 0.126 & 0.545 & 0.693\\
\hline
\end{tabular}
\captionof{table}{2D example ($\kappa = 1$).}

\end{minipage}


\begin{thebibliography}{1}
\bibitem{wachsmuth2016}
 D.\,Wachsmuth and F.\,P\"orner, arXiv:1603.05792 (2016)

\bibitem{fenics2016}
FEniCS, FEniCS project, http://www.fenicsproject.org/, 2016

\bibitem{beuchler2010}
 S.\,Beuchler, C.\,Pechstein and D.\,Wachsmuth, Comput Optim Appl. \textbf{51}, 883-908 (2010)
 
\bibitem{troeltzsch2010}
 F.\,Tr\"oltzsch, Optimal Control of Partial Differential Equations (American Mathematical Society, Providence, RI, 2010)
\end{thebibliography}
\end{document}